\documentstyle[aps,psfig]{revtex}
\input defDiag.sty

\begin{document}
\draft{
\title{Spectrum of stochastic evolution operators: polynomial basis
approach}

\author{C.P.Dettmann and Niels S\o ndergaard}
\address{Northwestern University, Department of Physics \&\ Astronomy\\
2145 Sheridan Road, Evanston, Illinois 60208 USA}
\author{Gergely Palla and G\'abor Vattay}
\address{Departement of Physics of Complex Systems, E\"otv\"os University\\
P\'azm\'any P\'eter s\'etany 1/A,
H-1117 Budapest, Hungary   }

\date{\today}

\maketitle

\begin{abstract}
The spectrum of the \evOper\  associated with a nonlinear stochastic flow
with additive noise is evaluated by diagonalization in
a polynomial basis. The method works for arbitrary noise strength.
In the weak noise limit we formulate a
new perturbative expansion  for the spectrum of
the stochastic \evOper\ in terms of expansions
around the classical periodic orbits.
The diagonalization of such operators is easier to implement than
 the standard Feynman diagram perturbation theory.
 The result is a stochastic analog of the Gutzwiller semiclassical
 spectral determinant
 with  the ``$\hbar$'' corrections computed to at least two orders more than
 what has so far been attainable in stochastic and quantum-mechanical
 applications, supplemented by the estimate for the late terms in the
 asymptotic saddlepoint expansions.
 \end{abstract}
 }

\section{Introduction}

The periodic orbit theory relates the spectrum of the Fokker-Plank operator and
its weighted \evOper\ generalizations to the periodic orbits
via trace formulas, \dzeta s and \fd s\rf{PG97,QCcourse}. 
For quantum mechanics the periodic orbit theory is exact on the
semiclassical level\rf{gutbook}, whereas the quintessentially quantum
effects such as creeping, tunneling
and diffraction have to be included as corrections. In particular,
the higher order $\hbar$ corrections can be computed
perturbatively by means of Feynman diagrammatic expansions\rf{alonso1}.

The \evOper\ formalism allows us to calculate long time averages
in a chaotic system in terms of the eigenvalues of \evOper s.
The simplest example is provided by the {\FPoper} 
\[
\Lop \rho(x')=\int dx\,\delta(f(x)-x')\rho(x)
\]
for a {\em deterministic} map $f(x)$ which maps a density distribution
$\rho(x)$ forward in time.
Our purpose here is to develop effective methods for computation of
spectra of {\em stochastic} \evOper s.
In case at hand, already a discrete time 1-dimensional 
discrete Langevin equation\rf{vk,LM94},
\begin{equation}
x_{n+1}=f(x_n)+\sigma\xi_n
\,,\label{Langevin}
\end{equation}
with $\xi_n$ independent normalized
random variables,
suffices to reveal the structure of the dependence on the noise.

We treat a chaotic system with 
external noise by 
replacing the deterministic evoluton $\delta$-function kernel 
by  the Fokker-Planck
kernel corresponding to (\ref{Langevin}),
a sharply peaked noise distribution function 
\beq
\Lop(x',x) =\delta_\sigma(f(x)-x')
\,.
\ee{Lnoise}
In \refref{noisy_Fred} we have treated the problem
of computing the spectrum of this operator by standard
field-theoretic Feynman diagram expansions; in
\refref{conjug_Fred} we offered a more elegant formulation
of the perturbation theory in terms of smooth conjugacies.
This time we 
evaluate the \evOper\ in an explicit polynomial 
basis. The procedure, which is relatively
easy to automatize, enables us to 
go three orders further in the perturbation theory;
in the language of Feynman diagrams of \refref{noisy_Fred},
in the new approach developed below we are
able to compute the perturbative corrections up 5-loop level,
as well as study the asymptotics of late terms in the
perturbative expansion.

The paper is organized as follows:
in refsect{FLOWS} we review the 
\evOper\ formalism for smooth flows.
In \refsect{ENTIRE} we explain the theorems that guarantee
that \Fd s for Axiom~$A$ systems are entire.
In refsect{NUMERICAL} we discuss the numerical tests of our
perturbative expansions.

\section{Matrix representation of {\FPoper}}
\label{ENTIRE}

%
We shall sketch here the basic ideas behind the
proofs that the \cFd s are entire,
without burdening the reader with too many technical details
(rigorous treatment is given in refs.~\cite{Ruelle76,frie,Rugh92}).
The main point
is that the \Fd s are entire
functions in any dimension, provided that

\noindent
1. the \evOper\ is {\em multiplicative} along the flow,
\\
\noindent
2. the symbolic dynamics is a {\em finite subshift}, 
\\
3. all cycle eigenvalues are {\em hyperbolic} (sufficiently bo\-und\-ed 
away from 1),
\\
4. the map (or the flow) is {\em real analytic}, {\em ie.}
   it has a piecewise  analytic continuation to a complex extension
   of the phase space.

As in physical applications one studies smooth dynamical observables,
we restrict the space that $\Lop$ acts on to smooth functions.
In practice ``real analytic'' means that 
all expansions are polynomial expansions.
In order to illustrate how this works in practice, we first work out
a simple example.

\subsection{Expanding maps with a single fixed point}

We start with the trivial
example of a repeller with only one expanding linear branch\rf{QCcourse}
\[  f(x) = \ExpaEig x \ \ \ \ |\ExpaEig| > 1 \ .\]
The action of the associated deterministic, noiseless \FPoper\ is 
\[
\Lop \phi (y) =  \int dx \delta(y-\ExpaEig x)\phi(x)
              = {1 \over |\ExpaEig|} \phi (y/\ExpaEig) \,. 
\] 
From this one immediately identifies the eigenfunctions and eigenvalues:
\beq
 \Lop \, y^n = {1 \over |\ExpaEig| \ExpaEig^{n}} y^n  \,,
\ \ \ \ n = 0,1,2,\ldots 
\ee{FP_eigs}
We note that the eigenvalues $\ExpaEig^{-n-1}$ 
fall off exponentially with $n$, and
that the trace of $\Lop$ is given by
\beq \tr \Lop = {1 \over |\ExpaEig|} \sum_{n=0}^\infty \ExpaEig^{-n} 
                  = {1 \over |\ExpaEig| (1-\ExpaEig^{-1})}
                   = {1 \over |f' - 1|} 
\,.
\ee{diagTr}
A similar result is easily obtained
for powers of $\Lop$, and for the \Fd\ one obtains:
\[
 \det(1-z\Lop) = \prod_{k = 0}^\infty \left(1-{z\over|\ExpaEig| \ExpaEig^{k}}
				      \right) =
        \sum_{k = 0}^\infty Q_k t^k 
\,,
\]
$ t =-z/|\ExpaEig|$, 
where the cumulants $Q_k$ are given explicitly by the Euler 
formula\rf{eule}
\beq
 Q_k = \frac{        1      }{1-\ExpaEig^{-1}} \;
         \frac{\ExpaEig^{-1}}{1-\ExpaEig^{-2}} \; \cdots \;
         \frac{\ExpaEig^{-k+1}}{1-\ExpaEig^{-k}}  \ \ \ .
\ee{Euler}
These coefficients decay asymptotically {\em faster} than exponentially, 
as $\ExpaEig^{-k(k-1)/2}$.
An intuitive way of comprehending the spectrum is to view the composition
with an expanding map as a ``smoothing" operator;
fast variations in the initial density distribution $\phi(x)$, corresponding to
high powers of $x$, are wiped out quickly by the smoothing operator.

In a suitable polynomial basis $\phi_n(x)$ the operator
has an explicit matrix representation
\[ (\Lop \phi)_n(x) = \sum_{m = 0}^\infty L_{nm}\phi_m(x) 
\,.
\]
In the single fixed-point example \refeq{FP_eigs}, 
$\phi_n = y^n$, and $ \Lop $ is diagonal, 
$\Lop_{nm}= \delta_{mn} \ExpaEig^{-n}/|\ExpaEig|$.
In general case, a matrix representation can be constructed 
by means of Cauchy complex contour integrals.

The simplest example of how the Cauchy formula
is employed is provided by a nonlinear inverse map 
$\psi=f^{-1}$,
$s = \sign{ \psi'}$
\[ \Lop \phi(w) 
  = \int d\!x \; \delta(w-f(x)) \phi(x) 
  = s \; \psi'(w)\; \phi(\psi(w)) 
\,.\]
Assume that $\psi$ is a contraction of the unit disk, i.e.
\[ |\psi(w)| < \theta < 1 \ \ \ \mbox{and} \ \ \ 
   |\psi'(w)| < C < \infty
  \  \ \ \mbox{for} \ \ \ |w|<1 
\,,
\]
and expand $\phi$ in a polynomial basis 
by means of the Cauchy formula
\beq
 \phi(x) = \sum_{n \geq 0} x^n \phi_n = \oint \frac{dw}{2 \pi i} \;\frac{\phi(w)}{w-x} 
\,,\quad
\phi_n = \oint \frac{dw}{2 \pi i}\; \frac{\phi(w)}{w^{n+1}} 
\,.
\ee{CauchL}

In this basis, $\Lop$ is a represented by the matrix
\beq
 \Lop \phi(w) = \sum_{m,n} w^m L_{mn} \phi_n 
\,,\quad
L_{mn} = \oint \frac{dw}{2 \pi i} \;\frac{s \; \psi'(w) (\psi(w))^n}{w^{m+1}}
\,.
\ee{diagLmat}
Taking the trace and summing we get:
\[ \tr \Lop = \sum_{n\geq 0} L_{nn} = \oint \frac{dw}{2 \pi i}\;
    \frac{s \; \psi'(w)}{w-\psi(w)} \]
This integral has but one simple pole at the unique fix point
$w^* = \psi(w^*) = f(w^*)$. Hence
\[ \tr \Lop = \frac{s \; \psi'(w^*)}{1-\psi'(w^*)} =
    \frac{1}{|f'(w^*)-1|} \]
in agreement with \refeq{diagTr}.

In practice we do not evaluate the matrix elements
\refeq{diagLmat}
by Cauchy integrals; instead, 
we use the polynomial basis\rf{Rugh92,CRR93}, 
with the left basis given by derivatives $\frac{z^m}{(m!)}$,
and the right basis by polynomials $\frac{\pde^m}{\pde x^m}$.
An analytic function can be written as
$F(y)=\sum F_m\frac{y^m}{m!}$, and the coeficients $F_m$
can be obtained by
$F_m=\left.\frac{\pde^m F(y)}{\pde y^m}\right|_{y=0}$, or in other
words, the dual basis of $\frac{y^m}{m!}$
is $\frac{\pde}{\pde^m}$. So the $L_{lk}$ matrix
element is obtained by acting the operator on $\frac{y^k}{k!}$
and differentiating the result $l$ times at $y=0$:
\bea
L_{m m'}^i &=& 
\left<\frac{\pde^m}{\pde y'^m}\right|\Lop_\Ssym{i}(y',y)
				\left|\frac{y^{m'}}{m'!}\right>
	\continue
	&=&
\left<\frac{\pde^m}{\pde y'^m}\right|\delta(y'+x_{i+1}-f(y+x_i))\left|\frac{y^{m'}}{m'!}\right>
\label{diagBtoL}
\eea

While it is not at all obvious that what is true for a single
fixed point should also apply to a Cantor set of periodic points,
the same asymptotic decay of expansion coefficients
is obtained when several expanding branches are involved. 
For a two-branch repeller, the procedure is the same, except
the integral \refeq{CauchL}
picks up a contribution from each branch.

From bounds on the elements $\Lop_{mn}$ one 
verifies\rf{Ruelle76,frie,Rugh92}
that they again fall off as $\ExpaEig^{-k^2/2}$,
concluding that the $\Lop$ eigenvalues fall off exponentially
for a general Axiom~$A$ one-dimensional map.

\section{\EvOper\ composition for piecewise-analytic expanding maps}

Suppose we can decompose the \evOper\ in
a sum of two operators $ \Lop =\Lop_0 + \Lop_1$.
Here we shall study a noisy map with a binary Markov partition
$\pS= \{\pS_0,\pS_1\}$.
We distinguish two branches of the map :
$f(x)=f_0(x)$ if $x \in \pS_0$,
and $f(x)=f_1(x)$ if $x \in \pS_1$.
In the case of Gaussian noise the corresponding operators are
\bea
\Lop_0(x',x)
	&=&
	\frac{1}{\sqrt{2\pi}\sigma}e^{-\frac{1}{2\sigma^2}(x'-f_0(x))^2}
\,, \qquad	x \in \pS_0
	\continue
\Lop_1(x',x)
	&=&
\frac{1}{\sqrt{2\pi}\sigma}e^{-\frac{1}{2\sigma^2}(x'-f_1(x))^2}
\,, \qquad	x \in \pS_1
\,.
\eea
If you visualize $\Lop$ as an [$2n\times 2n$] matrix,
in the weak noise limit 
the matrix elements for whom $x$ is in the other partition
are negligeable and can be set to zero,
so $\Lop_0$ and $\Lop_1$ are [$n\times 2n$] matrices.

For piecewise-analytic maps the decomposition is exact;
for stochastic operators is assumes that the overlap of the
two kernels is insignificant and can be neglected.
Then we can write:
\bea
\ln\det(1-z(\Lop_0+\Lop_1))
	&=& 
	\tr\ln(1-z(\Lop_0+\Lop_1))
	\continue
	&=& 
	-\sum_{n=1}^\infty\frac{z^n}{n}\tr(\Lop_0+\Lop_1)^n
\,.
\label{trl}
\eea
A contribution to this sum corresponds either to
a prime cycle $p = \Ssym{1} \Ssym{2} \cdots \Ssym{\cl{p}}$,
$\Ssym{i} \in \{0,1\}$, or to its repeat
\[
\tr {\Lop}_{p^r}= 
\tr( \Lop_\Ssym{1}\Lop_\Ssym{2}\cdots\Lop_\Ssym{\cl{p}})^r
\,.
\]
Expanding the $n$th power we get contributions from all 
symbol sequences of length $n$:
\bea
\tr(\Lop_0+\Lop_1)^n
= \sum_p \cl{p} \sum_{r=1}^\infty 
\delta_{n,\cl{p} r}
\tr (\Lop_p)^r,
\label{tracenp}
\eea
where $p$ denotes a prime cycle itinerary composed of $0$'s
and $1$'s, $\cl{p}$ is the prime cycle length and $r$ is the repetition
number. 
For instance,  $p=011$ then $\Lop_{011}=\Lop_0\Lop_1\Lop_1=\Lop_0\Lop_1^2$.
The cyclic property of trace yields the $\cl{p}$ factor.
Substituting into \refeq{trl} we obtain
\beq
\ln\det(1-z(\Lop_0+\Lop_1))
	=-\sum_p \sum_{r=1}^\infty \frac{z^{\cl{p}r}}{r}\tr\Lop_p^r
\,.
\ee{diagLnDet}
The sum over repeats yields a factorized formula of the \dzeta\ type
\beq
\det(1-z(\Lop_0+\Lop_1))
	=\prod_p\det(1-z^{\cl{p}}\Lop_p)
\,.
\ee{prod}
The operators are defined on piecewise monotonic maps, so there
is only one periodic orbit for a given prime cycle itinerary. 

\section{Saddle point expansions in terms of prime cycles}

In the weak noise limit the kernel is sharply peaked, so it
makes sense to expand it
in terms of the Dirac delta function and
its derivatives:
\bea
	\delta_\sigma(y)
	&=& 
	\sum_{m=0}^{\infty} {a_m \sigma^m \over m!} \, \delta^{(m)}(y) 
	\continue
	&=&
	\delta(y) + 
	a_2 {\sigma^2 \over 2} \delta^{(2)}(y) +
	a_3 {\sigma^3 \over 6} \delta^{(3)}(y) + \dots
	\,.
\label{delSigExp}
\eea
where
\[
	\delta^{(k)}(y) = {\pde^k \over \pde y^k} \delta(y)
	\,,
\]
and the coefficients $a_m$ depend on the choice of the kernel.
We have omitted the $\delta^{(1)}(y)$ term in the above because 
in our applications we shall impose
the saddle-point condition, that is, 
we shift $f$ by a constant to ensure that the noise peak corresponds
to $y=0$, so $\delta_\sigma^{'}(0)=0$.
For example, if $\delta_\sigma(y)$ is a Gaussian kernel,
it can be expanded as
\bea
	\delta_\sigma(y)
	&=& 
	{1 \over \sqrt{2 \pi \sigma^2}} e^{-{y^2/2\sigma^2} }
	=
	\sum_{n=0}^{\infty}
	\frac{\sigma^{2n}}{n!2^n} \delta^{(2n)}(y)
	\continue
	&=& 
	\delta(y) + {\sigma^2 \over 2} \delta^{(2)}(y) 
	 + {\sigma^4 \over 8} \delta^{(4)}(y) + \cdots
	\,.
\label{delGaussExp}
\eea

\section{\EvOper s in a matrix representation}
In this section our goal is to calculate noise corrections to the
leading eigenvalue of the Perron-Frobenius operator. Expression ($\ref{trl}$)
 showes, that in order to do that first we have to calculate $\tr{\Lop}^n$, 
and from equation ($\ref{tracenp}$) we see that $\tr{\Lop}^n$ should
 be generated from traces of $\Lop$ on {\em periodic orbits}. This section
 is a breaf review of how this  was carried out in practise using 
 the matrix representation of the Perron-Frobenius operator. 
 
If the coordinates of a prime cycle are $x_1,...,x_{\cl{p}}$,
the operator of a periodic orbit segment is
\beq
{\Lop}_\Ssym{i}(y',y)=\sum_{m=0}^{\infty}\frac{a_m\sigma^{m}}{m!}
\delta^{m}(y'+x_{i+1}-
f(x_i+y))
\,,
\label{li}
\eeq

where $x_{\cl{p}+1}=x_1$ and $a_m$ are the moments of the noise.
The full contribution to the trace from this periodic orbit is:
\[
\tr {\Lop}_{p}=\int dy_1...dy_p 
\Lop_\Ssym{\cl{p}}(y_{1},y_\cl{p})
\cdots
\Lop_\Ssym{2}(y_{3},y_2)
\Lop_\Ssym{1}(y_{2},y_1)
\,.
\]
The calculation of this contribution can be computerized if we represent 
(\ref{li}) in a matrix form using the polynomial basis
\refeq{diagBtoL}.
The matrix elements of $\Lop_\Ssym{i}$ are:
\beq
L_{m m'}^i=\left<\frac{\pde^m}{\pde y'^m}\right|\Lop_\Ssym{i}(y',y)
\left|\frac{y^{m'}}{m'!}\right>=\sum_{n=max(m'-m,0)}^{\infty}
\frac{a_n\sigma^n}{n!}B_{m+n,m'}^{(i)},
\label{BtoL}
\eeq
where the $B$ matrix is the representation of the noiseless operator:
\beq
B_{k k'}^i=\left<\frac{\pde^k}{\pde y'^k}\right|\delta(y'+x_{i+1}-f(y+x_i))\left|\frac{y^{k'}}{k'!}\right>
\eeq
If the Dirac-delta in $B$ acts on $\frac{y^{k'}}{k'!}$ we get:
\bea
& &\int dy\delta(y'+x_{i+1}-f(x_i+y))\frac{y^{k'}}{k'!}=
\frac{\left(f^{-1}(x_{i+1}+y')-x_i\right)^{k'}}{k'!|f'\left(f^{-1}(x_{i+1}+y')
\right)|}
 \nnu  \\
& &=\frac{\sign{f'}}{(k'+1)!}\frac{d}{dy'}\left(f^{-1}(x_{i+1}+y')-x_i
\right)^{k'+1}, 
\eea
where $\sign{f'}$ is
 the sign of $f'(f^{-1}(x_{i+1}+y'))$.
There is no summation for the branches of the prepimages of $f$, since 
each branche takes us back to a differernt surrounding, so at each 
point the branch is choosen by the previous point of the orbit.
The $k,k'$ matrix element of $B$ is the $k$th derivative of this at
$y'=0$:
\bea
B_{k,k'}^i&=&\frac{\sign{f'}}{(k'+1)!}
\left.\frac{\pde^{k+1}}{\pde y'^{k+1}}
		F^i(y')^{k'+1}
	\right|_{y'=0},
\label{Bmatrix} \\
F^i(y')&=&f^{-1}(x_{i+1}+y')-x_i
\nnu
\eea
Let us write $F^i(y')$ in power series of $y'$.
If $k'>k$ then after the derivates are taken 
every term in $B_{k,k'}$ has at least a $(y')^{k'-k}$
coefficient. We have to evaluate the derivative at $y'=0$, so
 we see that $B$ is a triangular matrix,
$B_{k,k'}=0$ if $k'>k$. 

\subsection{Multinomials}

The non-zero matrix elements are\rf{Abramowitz} 
\beq
\left(\sum_{l=1}^{\infty}\frac{x_l}{l!}t^l\right)^m=m!\sum_{n=l}^{\infty}
\frac{t^n}{n!}\sum(n|a_1 a_2 ... a_n)'x_1^{a_1}x_2^{a_2}...x_n^{a_n}, 
\label{Abram}
\eeq
where the sum $\left(\sum\right)$ goes over all non-negative integers
such that:
\beq
a_1+2a_2+...+n a_n=n 
	\,,\qquad
a_1+a_2+...+a_n=m
\eeq
and the multinomial coefficient is:
\beq
(n|a_1 a_2 ... a_n)'=\frac{n!}{(1!)^{a_1}a_1!(2!)^{a_2}a_2!...(n!)^{a_n}a_n!}
\eeq
If we expand $F^i(y')$ in a Taylor series, the constant term is
zero because of $f^{-1}(x_{i+1})=x_i$. So we can write:
\beq
F^i(y')=\frac{y'}{\Lambda}+\sum_{l=2}^{\infty}\frac{F^i_l}{l!}y'^l
\eeq
We apply the formula ($\ref{Abram}$) to $F^i(y')$ with power $k'+1$:
\beq
(F^i(y'))^{k'+1}=(k'+1)!\sum_{n=k'+1}^{\infty}\frac{y'^n}{n!}\sum(n|a_1 a_2 ... a_n)'\frac{1}{\Lambda^{a_1}}(F^i_2)^{a_2}...(F^i_n)^{a_n}.
\eeq
In the $(k+1)$th derivative of this at $y'=0$ only the $n=k+1$ term is
 non-zero. So considering only one branch, with positve sign of $f'_b$:
\bea
& &B^i_{k,k'}=\sum(k+1|a_1 a_2 ... a_{k+1})'\frac{1}{\Lambda^{a_1}}(F^i_2)^{a_2}...
(F^i_{k+1})^{a_{k+1}} \\
& &a_1+2a_2+...(k+1)a_{k+1}=k+1 \label{a1} \\
& &a_1+a_2+...a_{k+1}=k'+1 \label{a2} \\
& &(k+1|a_1 a_2 ... a_{k+1})'=\frac{(k+1)!}{(1!)^{a_1}a_1!(2!)^{a_2}a_2!...((k+1)!)^{a_{k+1}}a_{k+1}!} \label{coef}
\eea
(Here $F^i_l$ is the $l$th derivative of $F^i$  on the particular
branch which we have to choose.)
For the diagonal and the nearest off-diagonals the matrix elements read as:
\bea
& &B_{m m}=\frac{1}{\Lambda^{m+1}}, \mbox{\hspace{0.5cm}}m=0,1,... \\
& &B_{m+1,m}=\frac{(m+2)(m+1)}{2}\frac{F^i_2}{\Lambda^m} \\
& &B_{m+2,m}=
	\frac{(m+3)!}{8(m-1)!}\frac{(F^i_2)^2}{\Lambda^{m-1}}+ \frac{(m+3)!}
{6(m+1)!}\frac{F^i_3}{\Lambda^m}\\
& &B_{m+3,m}=\frac{(m+4)!}{48(m-2)!}\frac{(F_2^i)^3}{\Lambda^{m-2}}+\frac{(m+4)!}{12(m-1)!}\frac{F^i_2F^i_3}{\Lambda^{m-1}}+\frac{(m+4)!}{24m!}\frac{F^i_4}{\Lambda^m}\\
& &\mbox{etc...} \nnu
\,,
\eea

\subsection{Noiseless case}
In the noiseless case the $B^i$ matrices are the representations
 of the \FPoper. Since they are triangular, their eigenvalues are:
\beq
\lambda_m=B_{m m}=\frac{1}{\Lambda^{m+1}}\mbox{\hspace{0.5cm}}.
\eeq
Multiplied, the triangular matrices yield a
triangular matrix, and the diagonal elements of the result can be 
obtained simply by multiplying the two corresponding diagonal elements.
 This results that the trace of the $\Lop$ on a periodic orbit is
 the following:
\bea
\tr \Lop_p
	&=& \tr L_\Ssym{1}L_\Ssym{2}\cdots L_\Ssym{\cl{p}}
	    =\sum_{m=0}^{\infty} \frac{1}{|\Lambda_p|^{(m+1)}}
	\continue
	&=&
	\frac{1}{|\Lambda_p|}
	\sum_{m=0}^{\infty}\frac{1}{\Lambda_p ^m}
	=\frac{1}{|\Lambda_p|\left(1-1/\Lambda_p\right)}
	=\frac{1}{|1-\Lambda_p|}
\eea
The absolute value is a consequence of the $\sign{f'}$ 
factor in \refeq{Bmatrix}.
For the total trace, using ($\ref{tracenp}$) we get back the standard
deterministic trace formula\rf{QCcourse}
\beq
\tr{\Lop}^n=
\sum_p \cl{p} \sum_{r=1}^\infty
\delta_{n,\cl{p} r}
\frac{1}{|1-\Lambda_p^r|}
\eeq

\subsection{Numerical tests}

Here we continue the calculations of Sect.~? of \refref{conjug_Fred},
where more details and discussion may be found.  
We test our perturbative expansion on the repeller
of the 1-dimensional map

\beq 
f(x)=20\left(\frac{1}{16}-\left(\frac{1}{2}-x\right)^4\right).
\ee{testQuartic}
This repeller is a nice example of an ``Axiom~$A$'' expanding
system of bounded nonlinearity and complete binary symbolic dynamics, 
for which the deterministic \evOper\ eigenvalues
converge super-exponentially with the cycle length.

The inverse of $f$  has two branches:
\beq
f^{-1}(y)=\left[\frac{1}{2}-\left(\frac{1}{16}-\frac{y}{20}\right)^{\frac{1}{4}},\frac{1}{2}
+\left(\frac{1}{16}-\frac{y}{20}\right)^{\frac{1}{4}}\right].
\eeq
Due to the symmetry of \refeq{testQuartic} around $x=0.5$,
the derivatives on the two inverse branches have the same absolute value,
but opposite signs. 
We compute the leading eigenvalue of the {\evOper}
(the repeller escape rate) in the presence of Gaussian noise, using three
complementary approaches.  The perturbative result in terms
of periodic orbits and the weak noise corrections
is compared to the eigenvalue computed by a numerical
lattice discretization in \refref{noisy_Fred}. In the preceding
paper (\refref{conjug_Fred}) we compared the numerical eigenvalue with the $\sigma^4$ result
and estimated the coefficient of $\sigma^6$ to be approximately 2700.  Here we
compute the order $\sigma^6$ coefficient $2076.47\ldots$
to 14 digits accuracy, as well as the $\sigma^8$ and $\sigma^{10}$
coefficients. Furthermore, we estimate the asymptotic form of the
$\sigma^{2m}$ coefficient.

The numerical calculations of $\tr \Lop^n$ proceeds as follows:
\begin{enumerate}
\item
	Creating periodic orbits up to length 10 using the method 
	of iterating backwards\rf{QCcourse}
\item{Computing [$24\times 24$] $B$ matrices at each point of orbits}
\item{From the $B$ matrices computing [$16\times 16$] $L$ matrices at each point of
 orbits, up to $\sigma^{10}$ corrections in each matrix-element}
\item{Multiplying the $L$ matrices. At shorter orbits repeating this $r$ 
 times.($r\cl{p}=n$)}
\item
Computing the traces of the result matrices, multiplying them by $\cl{p}$ 
and adding them up to get $\tr\Lop^n$ as in \refeq{tracenp}

\item Create cummulants from the traces.
\item Find zeros.
\end{enumerate}

At $n=2$ [$20\times 20$] sized $L$ (and corresponding to that
[$28\times 28$] sized $B$) matrices were used because of slower convergence in 
$\tr\Lop^n$ while enlarging the matrix size. At $n=1$ 
[$26\times 26$] is the size in $L$ one has to reach to get decent result for
$\tr\Lop$.

\subsection{From traces to $\sigma^{2m}$ corrections}
To get the leading eigenvalue of the Perron-Frobenius operator from the
traces of ${\cal L}^n$, we follow the method outlined in\rf{conjug_Fred}.
The trace of $\Lop^n$ for $n=1,2,...,10$ can be written as:
\bea
\tr\Lop^n=\sum_{j=0}^{\infty}C_{n j}\sigma^j,
\label{ce}
\eea
where we know the value of the first five non-zero 
$C_{n 0}, \cdots, C_{n,10}$ . 
The cumulants $Q_{n j}$ in
\bea
\det(1-z\Lop)=1-\sum_{n=1}^{\infty}\sum_{j=0}^{\infty}Q_{n j}z^n\sigma^j
\label{qum}
\eea
are obtained recursively  as
\bea
Q_{n m}&=&\frac{1}{n}\left(C_{n m}
	   - \sum_{k=1}^{n-1}\sum_{l=0}^{m}Q_{k,m-l}C_{n-k,l}
		     \right)
\,.
\eea
For Gaussian noise only even $l$ terms contribute.
Let $z_0$ be the solution of the noiseless condition
$\left.\det(1-z\Lop)\right|_{\sigma=0}=0$.
Expanding the 
spectral determinant around $z_0$ and $\sigma^2=0$, we write:
\bea
\det(1-(z+z_0)\Lop)&=&F-F_{10}z-F_{02}\sigma^2-F_{20}z^2-F_{12}z\sigma^2-F_{04}\sigma^4 \nnu \\ & &-{\cal O}(z^3,\sigma^6)
\label{sor}
\eea
where the coefficients can be obtained from the cumulants as:
\bea
\matrix{F&=&1-\sum_{m=1}^nQ_{m,0}z_0^m & & F_{10}&=&\sum_{m=1}^nmQ_{m,0}z_0^{m-1} \cr & & & & &
 & \cr F_{02}&=&\sum_{m=1}^nQ_{m,2}z_0^m & & F_{20}&=&\frac{1}{2}\sum_{m=2}^nm(m-1)Q_{m,0}z_0^{m-2} \cr & & & & & & \cr F_{12}&=&\sum_{m=1}^nmQ_{m,2}z_0^{m-1} & & F_{04}&=&
\sum_{m=1}^nQ_{m,4}z_0^m}
\eea
If we want to calculate only the $\sigma$ corrections just up to fourth order,
first we have to solve  this equation:
\bea
& &F-F_{10}(z_2\sigma^2+z_4\sigma^4)-F_{02}\sigma^2-F_{20}(z_2\sigma^2+
z_4\sigma^4)^2\nnu \\ & &-F_{12}(z_2\sigma^2+z_4\sigma^4)\sigma^2-F_{04}\sigma^4=0.
\eea
Since $\sigma$ is a parameter, the coefficients of $\sigma^0$, or
 $\sigma^2$ or $\sigma^4$ have to vanish:
\bea
 F&=&0 \nnu \\ -F_{10}z_2-F_{02}&=&0 \nnu \\-F_{10}z_4-F_{20}z_2^2-F_{12}z_2
-F_{04}&=&0 \nnu 
\eea 
The first equation gives solution for $z_0$, which numerically was
done by Newton's method. From the second we can
express $z_2$, from the third $z_4$:
\bea
z_2=-\frac{F_{02}}{F_{10}}\mbox{\hspace{1cm}}z_4=-\frac{F_{20}z_2^2+F_{12}z_2+F_{04}}{F_{10}}
\label{ifi}
\eea
To get the eigenvalue we have to reciprocate $z$:
\bea
\nu_0+\nu_2\sigma^2+\nu_4\sigma^4=\frac{1}{z_0+z_2\sigma^2+z_4\sigma^4}&=&
 \frac{1}{z_0}-\frac{z_2}{z_0^2}\sigma^2-\left(\frac{z_4}{z_0^2}-\frac{z_2^2}{z_0^3}
\right)\sigma^4 \nnu \\ & &+{\cal O}(\sigma^6)
\eea
We see that:
\bea
\nu_0=\frac{1}{z_0},\mbox{\hspace{0.5cm}}\nu_2=-\frac{z_2}{z_0^2},
\mbox{\hspace{0.5cm}}\nu_4=-\frac{z_4}{z_0^2}+\frac{z_2^2}{z_0^3},
\eea
or by substituing into these \refeq{ifi}:
\bea
\nu_0 &=&\frac{1}{z_0}
	\,,\qquad 
\nu_2=\frac{F_{02}}{F_{10}}\nu_0^2,
	\continue
\nu_4 &=&
	\frac{1}{F_{10}^3} (F_{20}F_{02}^2-F_{12}F_{02}F_{10}+
	F_{04}F_{10}^2+F_{02}^2F_{10}\nu_0)
	\nu_0^2.
\eea
To go up further in the powers of $\sigma$ corrections, we have to
 do the same as we have done for order four, the details of this
 is in the appendix.

\begin{table}
\begin{tabular}{clll}\hline
$n$&$\nu_0$&$\nu_2$&$\nu_4$\\\hline
1&0.308&0.42&~2.2\\
2&0.37140&1.422&32.97\\
3&0.3711096&1.43555&36.326\\
4&0.371110995255&1.435811262&36.3583777\\
5&0.371110995234863&1.43581124819737&36.35837123374\\
6&0.371110995234863&1.43581124819749&36.358371233836\\\hline
$n$&$\nu_6$&$\nu_8$&$ $\\\hline
1 & 17.4             & 168.0            & \\
2&  1573.3           & 112699.9         & \\
3&   2072.9           &  189029.0         & \\
4 & 2076.479         &  189298.8         & \\
5 & 2076.4770492     &  189298.12802     & \\
6 & 2076.47704933320 &  189298.128042526 & \\
7 & 2076.47704933321 &  189298.128042526 & \\
\hline
\end{tabular}
\caption{Significant digits of the leading deterministic eigenvalue
and its $\sigma^2, \cdots,  \sigma^{8}$ coefficients, calculated from the
spectral determinant as a function of the cycle truncation length $n$.
Note the super-exponential convergence of all coefficients. 
We have computed all cycles to length 10, but contributions
of those longer than $n=6$ lie below the machine precision.
\label{tabPert}}
\end{table}

The perturbative corrections to the leading eigenvalue
(escape rate) of the weak-noise \evOper\ are given in \reftab{tabPert},
showing super-exponential convergence with the truncation 
cycle length $n$.  The super-exponential convergence
has been proven for the deterministic,
$\nu_0$ part of the eigenvalue\rf{grothi,Rugh92}, but has not been
studied for noisy kernels.
It is seen that a good first approximation
is obtained already at $n=2$, using only 3 prime cycles,
and $n=6$ (23 prime cycles in all) is in this example
sufficient to exhaust the limits of double precision arithmetic.
The exact value of $\nu_6 = 2076.47\ldots$ is not wildly different
to  our previous numerical estimate\rf{conjug_Fred}  of 2700.

\section{Testing the results}
\begin{figure}[hbt]
\centerline{\psfig{file=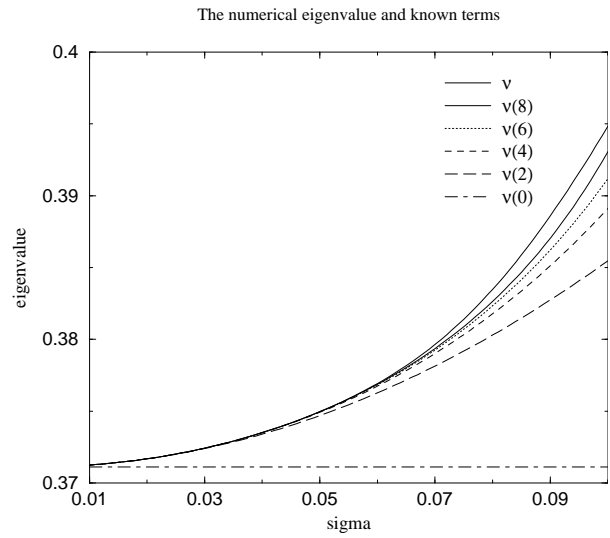,width=8cm}}
\caption{
The numerical eigenvalue and the known terms:$\nu(n)=\sum_{k=0}^{n/2}
\nu_{2k}\sigma^{2k}$\label{nufig1}}
\end{figure}

\begin{figure}[hbt]
\centerline{\psfig{figure=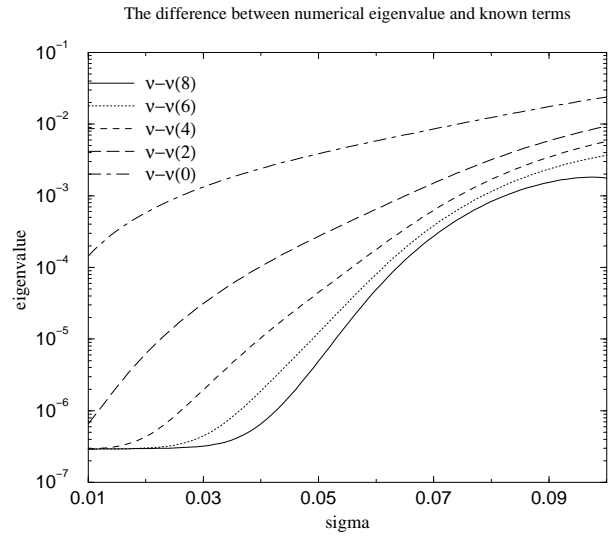,width=8cm}}
\caption{
The difference between the numerical eigenvalue and the known terms.
\label{nufig2}}
\end{figure}

\begin{figure}[hbt]
\centerline{\psfig{figure=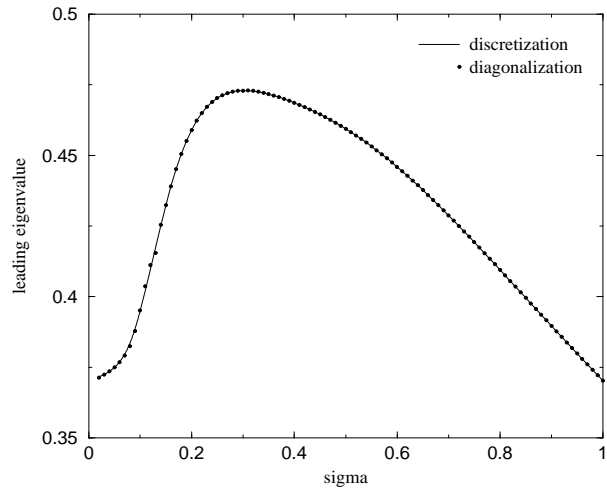,width=8cm}}
\caption{
The result of numerical integration and the eigenvalue fit
 together well in the large $\sigma$ region.
\label{dida}}
\end{figure}

In previous papers the eigenvalue was computed numerically. 
In \reffig{nufig1} to \reffig{nufig2} we check how well
our results fit that.
The present result fits the data of the numerical
 discretisation method well in the small noise region, but not 
 above $sigma=0.1$. At given sigma, a matrix representation 
of the Perron-Frobenius operator can be obtained numerically, 
at each matrix element carrying out a  numerical integration:
\bea
L_{l,k}=\left.\frac{\partial^l}{\partial y^l}\left[\frac{1}{\sqrt{2\pi}\sigma k!}\int dze^{-\frac{(y-f(z))^2}{2{\sigma}^2}}z^k\right]\right|_{y=0}.
\eea
The leading eigenvalue of this matrix fits well the discretisation data
well at large sigma, as it is shown on \reffig{dida}. 
\begin{figure}[hbt]
\centerline{\psfig{figure=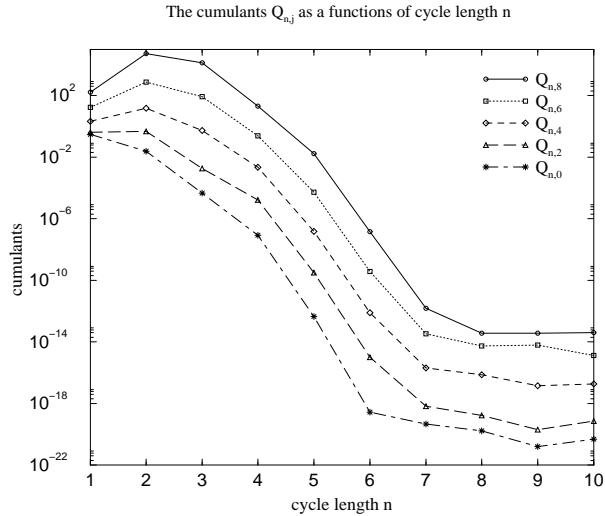,width=8cm}}
\caption{ The generalized cumulants $Q_{n j}$ as a function
of cycle length $n$ for $j=0,2,4,6,8$. $Q_{n 0}$ is
the cumulant of the noiseless case. Superexponential
convergence can be observed until cycle length $n=6$,
then numerical errors take over.\label{qumfig}}
\end{figure}

\section{Summary and outlook}

We have formulated a perturbation theory of
stochastic trace formulas based on a polynomial basis
matrix representations, expanded 
around infinitely many chaotic saddle points 
(unstable periodic orbits).

We note in passing that for
1-$d$ repellers a diagonalization of an explicit
truncated $\Lop_{mn}$ matrix yields many more eigenvalues than
the cycle expansions\rf{CCR,Rugh92}. The reasons why one persists
anyway in using the periodic orbit theory are partially aestethic,
and partially pragmatic. Explicit $\Lop_{mn}$ demands explicit
choice of a basis and is thus non-invariant, in contrast to
cycle expansions which utilize only the invariant information
about the flow. In addition, we do not know how to construct
$\Lop_{mn}$ for a realistic flow, such as the 3-disk problem, 
while the periodic orbit formulas are general and
straightforward to apply.

\section{Acknowledgements}
G.V. thanks the Hungarian Ministry of Education, OMFB, OTKA T25866/F17166
and the Humboldt foundation the financial support. G. V. thanks
Bruno Eckhardt  the  cordial hospitalty at the Department of Physics of
the Philipps-Universit\"at Marburg where some of the results
have been derived. G.P. thanks the support of FKFP0159/1997.


\appendix
\section{Algebra}
To go up further in the powers of $\sigma$ corrections in $\nu$, we have to
 do the same as we have done for order four. First expand ($\ref{sor}$):
\bea
\det(1-(z+z_0)\Lop)
&=&F-F_{10}z-F_{02}\sigma^2-F_{20}z^2-F_{12}z\sigma^2-F_{04}\sigma^4 \nnu \\
& &-F_{30}z^3-F_{22}z^2\sigma^2-F_{14}z\sigma^4-F_{06}\sigma^6 \nnu \\
& &-F_{40}z^4-F_{32}z^3\sigma^2-F_{24}z^2\sigma^4-F_{16}z\sigma^6-F_{08}\sigma^8\nnu \\
\label{alap}
\eea
The coefficients are:
\bea
\matrix{ F_{30}&=&\frac{1}{6}\sum_{m=3}^n\frac{m!}{(m-3)!}Q_{m,0}z_0^{m-3}& &
F_{22}&=&\frac{1}{2}\sum_{m=2}^n\frac{m!}{(m-2)!}Q_{m,2}z_0^{m-2} \cr & & & & & & \cr
F_{14}&=&\sum_{m=1}^nmQ_{m,4}z_0^{m-1} & & F_{06}&=&\sum_{m=1}^nQ_{m,6}z_0^m
\cr & & & & & & \cr F_{40}&=&\frac{1}{24}\sum_{m=4}^n\frac{m!}{(m-4)!}Q_{m,0}z_0^{m-4} & &
F_{32}&=&\frac{1}{6}\sum_{m=3}^n\frac{m!}{(m-3)!}Q_{m,2}z_0^{m-3} \cr & & & & & &
 \cr  
F_{24}&=&\frac{1}{2}\sum_{m=2}^n\frac{m!}{(m-2)!}Q_{m,4}z_0^{m-2} & &
F_{16}&=&\sum_{m=1}^nmQ_{m,6}z_0^{m-1} \cr & & & & & & \cr
 F_{08}&=&\sum_{m=1}^nQ_{m,8}z_0^m & & & & }
\eea
The equation to be solved is:
\bea
& &F-F_{10}(z_2\sigma^2+z_4\sigma^4+z_6\sigma^6+z_8\sigma^8)-
F_{02}\sigma^2-F_{20}(
z_2\sigma^2+z_4\sigma^4+z_6\sigma^6+z_8\sigma^8)^2\nnu \\& &
-F_{12}(z_2\sigma^2+z_4\sigma^4+z_6\sigma^6+z_8\sigma^8)\sigma^2
-F_{04}\sigma^4-F_{30}(z_2\sigma^2+z_4\sigma^4+z_6\sigma^6+z_8\sigma^8)^3 
\nnu \\& &-F_{22}(z_2\sigma^2+z_4\sigma^4+z_6\sigma^6+z_8\sigma^8)^2\sigma^2
-F_{14}(z_2\sigma^2+z_4\sigma^4+z_6\sigma^6+z_8\sigma^8)\sigma^4
-F_{06}\sigma^6 \nnu \\& & 
-F_{40}(z_2\sigma^2+z_4\sigma^4+z_6\sigma^6+z_8\sigma^8)^4
-F_{32}(z_2\sigma^2+z_4\sigma^4+z_6\sigma^6+z_8\sigma^8)^3\sigma^2 \nnu \\ & &
-F_{24}(z_2\sigma^2+z_4\sigma^4+z_6\sigma^6+z_8\sigma^8)^2\sigma^4 
-F_{16}(z_2\sigma^2+z_4\sigma^4+z_6\sigma^6+z_8\sigma^8)\sigma^6\nnu \\ & &
-F_{08}\sigma^8=0
\eea
The solutions for $z_6$ and $z_8$ are:
\bea
z_6&=&-\frac{1}{F_{10}}
(2F_{20}z_2z_4+F_{12}z_4+F_{30}z_2^3+F_{22}z_2^2+F_{14}z_2+F_{06}) \\
z_8&=&-\frac{1}{F_{10}}\left(F_{20}(2z_2z_6+z_4^2)+F_{12}z_6+3F_{30}z_2^2z_4
+2F_{22}z_2z_4\right. \nnu \\ & & \left.+F_{14}z_4+F_{40}z_2^4+F_{32}z_2^3+F_{24}z_2^2+F_{16}z_2+F_{08}\right)
\eea
The connection between the $\nu$'s and the $z$'s:
\bea
& &\nu_0+\nu_2\sigma^2+\nu_4\sigma^4+\nu_6\sigma^6+\nu_8\sigma^8=\frac{1}{
z_0+z_2\sigma^2+z_4\sigma^4+z_6\sigma^6+z_8\sigma^8} \nnu \\
&=&
=
\frac{1}{z_0}\left[1-\left(\frac{z_2}{z_0}\sigma^2+\frac{z_4}{z_0}\sigma^4+
\frac{z_6}{z_0}\sigma^6+\frac{z_8}{z_0}\sigma^8\right)\nnu \right.\\ & & \left.
+\left(\frac{z_2}{z_0}\sigma^2+\frac{z_4}{z_0}\sigma^4+
\frac{z_6}{z_0}\sigma^6+\frac{z_8}{z_0}\sigma^8\right)^2 -\left(\frac{z_2}{z_0}\sigma^2+\frac{z_4}{z_0}\sigma^4+
\frac{z_6}{z_0}\sigma^6+\frac{z_8}{z_0}\sigma^8\right)^3 \right. \nnu \\ & &\left.+\left(\frac{z_2}{z_0}\sigma^2+\frac{z_4}{z_0}\sigma^4+
\frac{z_6}{z_0}\sigma^6+\frac{z_8}{z_0}\sigma^8\right)^4...\right] \nnu \\
&=&\frac{1}{z_0}-\frac{z_2}{z_0^2}\sigma^2-\left(\frac{z_4}{z_0^2}-\frac{z_2^2}{z_0^3}
\right)\sigma^4-\left(\frac{z_6}{z_0^2}-2\frac{z_2z_4}{z_0^3}+\frac{z_2^3}{z_0^4}\right)\sigma^6 \nnu \\ & &-\left(\frac{z_8}{z_0^2}-2\frac{z_2z_6}{z_0^3}-\frac{z_4^2}{z_0^3}+
3\frac{z_2^2z_4}{z_0^4}-\frac{z_2^4}{z_0^5}\right)\sigma^8 
\eea
We can say that:
\bea
\nu_6&=&-\frac{z_6}{z_0^2}+2\frac{z_2z_4}{z_0^3}-\frac{z_2^3}{z_0^4} \nnu \\
&=&\left(2F_{02}^3F_{20}^2-3F_{02}F_{10}F_{12}F_{20}+2F_{02}
F_{04}F_{10}^2F_{20}+F_{02}F_{10}^2F_{12}^2\right.
\nnu \\ & &-F_{04}F_{10}^3F_{12}-F_{30}F_{02}^3
F_{10}+F_{22}F_{02}^2F_{10}^2-F_{02}F_{10}^3F_{14}
+F_{06}F_{10}^4\nnu \\ & & \left.+2(F_{02}^3F_
{10}F_{20}-F_{02}^2F_{10}^2F_{12}+F_{02}F_{04}F_{10}^3)\nu_0+F_{02}^3F_{10}^2\nu_0^2\right)\frac{\nu_0^2}{F_{10}^5} \\
\nu_8&=&-\frac{z_8}{z_0^2}+2\frac{z_2z_6}{z_0^3}+\frac{z_4^2}{z_0^3}-3\frac{z_2^2z_4}{
z_0^4}+\frac{z_2^4}{z_0^5} \nnu \\
&=&\left[F_{10}^4\left(F_{08}F_{10}^2-F_{06}F_{10}F_{12}+
F_{04}F_{12}^2-F_{04}F_{10}F_{14}+F_{04}^2F_{20}\right.\right. \nnu \\
& &\mbox{\hspace{1cm}}\left.+F_{04}^2F_{10}\nu_0\right) \nnu \\
& &+F_{02}^4\left(F_{04}F_{10}^2-5F_{20}^3-5F_{10}F_{20}F_{30}+
5F_{10}F_{20}^2\nu_0-2F_{10}^2F_{30}\nu_0 \right. \nnu \\
& &\mbox{\hspace{1cm}}\left.+3F_{10}^2F_{20}\nu_0^2+F_{10}^3\nu_0
\right) \nnu \\
& &+F_{02}^2F_{10}^2\left(-3F_{10}F_{14}F_{20}+6F_{04}F_{20}^2-3F_{10}F_{12}F_{22}+F_{10}^2F_{24} \right.\nnu \\ & & \mbox{\hspace{1.5cm}}
-3F_{04}F_{10}F_{30}-2F_{10}^2F_{14}\nu_0+
6F_{04}F_{10}F_{20}\nu_0 \nnu \\ & &\mbox{\hspace{1.5cm}}\left.
+3F_{04}F_{10}^2\nu_0^2+6F_{12}^2F_{20}+
3F_{12}^2F_{10}\nu_0\right)\nnu \\
& & -F_{02}F_{10}^3\left(F_{13}^3+F_{12}\left(-2F_{10}F_{14}+6F_{04}F_{20}+
4F_{04}F_{10}\nu_0\right)\right. \nnu \\ & &\mbox{\hspace{1.5cm}}
\Bigl.+F_{10}\left(F_{10}F_{16}-2F_{06}F_{20}-2F_{04}F_{22}
-2F_{06}F_{10}\nu_0\right)\Bigr)\nnu \\ 
& &-F_{02}^3F_{10}\Bigl(F_{10}\left(-4F_{20}F_{22}+F_{10}F_{32}-
2F_{10}F_{22}\nu_0\right)\Bigr.\nnu \\
& &\mbox{\hspace{1.5cm}}\Bigl. \left. +F_{12}(10F_{20}^2-4F_{10}F_{30}+
8F_{10}F_{20}\nu_0+3F_{10}^2\nu_0^2)\right)\Bigr]\frac{\nu_0^2}{F_{10}^7}
\eea

\section*{References} 

\end{document}